\documentclass[runningheads]{llncs}
\usepackage{graphicx}
\usepackage[symbol]{footmisc}

\usepackage{float}
\usepackage{mathtools}
\usepackage{amssymb, amsmath, latexsym}
\usepackage{nicefrac}
\usepackage{hhline}
\usepackage{multirow}

\usepackage{pifont}
\usepackage[colorinlistoftodos,bordercolor=orange,backgroundcolor=orange!20,linecolor=orange,textsize=scriptsize]{todonotes}
\usepackage{algorithm}\usepackage{algpseudocode}

\graphicspath{ {images/} }
\def \R {\mathbb R}

\def\R{\mathbb{R}}

\def\R{\mathbb R}

\begin{document} 
\title{Near-Optimal Decentralized Algorithms for Saddle Point Problems over Time-Varying Networks\thanks{The research of A. Beznosikov, A. Rogozin and A. Gasnikov was supported by Russian Science Foundation (project No. 21-71-30005).}}
\titlerunning{Near-Optimal Algorithms for SPPs over Time-Varying Networks}
%
\author{Aleksandr Beznosikov\inst{1,2} \and
 Alexander Rogozin\inst{1,2} \and
 Dmitry Kovalev\inst{5} \and
Alexander Gasnikov\inst{1,3,4}
}
\authorrunning{A. Beznosikov, A. Rogozin, D. Kovalev, A. Gasnikov}
%
\institute{Moscow Institute of Physics and Technology, Dolgoprudny, Russia \and
Higher School of Economics, Russia  \and
 Institute for Information Transmission Problems RAS, Russia \and
Caucasus Mathematical Center, Adyghe State University, Russia \and
King Abdullah University of Science and Technology, Saudi Arabia
}
\maketitle              
\begin{abstract}

Decentralized optimization methods have been in the focus of optimization community due to their scalability, increasing popularity of parallel algorithms and many applications. In this work, we study saddle point problems of sum type, where the summands are held by separate computational entities connected by a network. The network topology may change from time to time, which models real-world network malfunctions. We obtain lower complexity bounds for algorithms in this setup and develop near-optimal methods which meet the lower bounds.

\end{abstract}
\section{Introduction}

Distributed algorithms are an important part of solving many applied optimization problems \cite{shalev2014understanding,mcdonald2010distributed,mcmahan2017communication}. They help to parallelize the computation process and make it faster. In this paper, we focus on the distributed methods for the saddle point problem:
\begin{equation}
    \label{SPP}
    \min_{x \in \mathcal{X}} \max_{y \in \mathcal{Y}} f(x,y) := \frac{1}{M} \sum\limits_{m=1}^{M} f_m(x,y).
\end{equation}
In this formulation of the problem, the original function $f$ is divided into $M$ parts, each of part $f_m$ is stored on its own local device. Therefore, only the device with the number $m$ knows information about $f_m$. Accordingly, in order to obtain complete information about the function $f$, it is necessary to establish a communication process between devices. This process can be organized in two ways: centralized and decentralized. In a centralized approach, communication takes place via a central server, i.e. all devices can send some information about their local $f_m$ function to the central server, the server collects information from the devices and does some additional calculations, and then can send new information or request to the devices. Then the process continues. With this approach, one can easily write centralized gradient descent for distributed sum minimization: $\min_{x} g(x) := \frac{1}{M} \sum_{m=1}^{M} g_m(x)$. All devices compute local gradients in the same current point and then send these gradients to the server, in turn, the server averages the gradients and makes a gradient descent step, thereby obtaining a new current point, which it sends to the devices. Centralized methods for \eqref{SPP} are discussed in detail, for example, in \cite{beznosikov2020local}.  However, centralized approach has several problems, e.g. synchronization drawback or high requirements to the server. Possible approach to deal with these drawbacks is to use decentralized architecture \cite{bertsekas1989parallel}.  In this case, there is no longer any server, and the devices are connected into a certain communication network and workers are able to communicate only with their neighbors and communications are simultaneous. The most popular and frequently used communication methods are the gossip protocol \cite{kempe2003gossip,boyd2006randomized,nedic2009distributed} and accelerated gossip protocol \cite{scaman2017optimal,ye2020multi}. In the gossip protocol, nodes iteratively exchange data with their immediate neighbors using a communication matrix and in this way the information diffuses over the network. Decentralized algorithms are already widely developed for minimization problems, but not for saddle point problems. Meanwhile, saddle-point problems have a lot of applied applications, including those that require distributed computing. These are the already well-known and classic matrix game and Nash equilibrium\cite{GT-book,facchinei2007finite}, as well as modern problems in adversarial training \cite{Arjovsky_et_al2017,goodfellow2014generative}, image deconvolution \cite{chambolle2011first} and reinforcement, statistical learning \cite{pmlr-v119-jin20f,Abadeh_et_al_2015}.

This paper closes some of the open questions in decentralized saddle point problems.


\subsection{Our contribution}

In particular, our contribution can be briefly described as follows

\textit{Lower bounds.} We present lower bounds for decentralized smooth strongly-convex-strongly-concave and convex-concave saddle-point problems on the time-varying networks. The lower bounds are derived under the assumption that the network is always a connected graph.

\textit{Near-Optimal algorithm.} The paper constructs a near-optimal algorithm that meets the lower bounds. The analysis of the algorithm is carried out for smooth strongly-convex-strongly-concave and convex-concave saddle-point problems

See our results in the column "time-varying" of Table 1.

\subsection{Related works} 

Our work is one of the first dedicated to decentralized saddle problems over time-varying networks. Among other works, we can highlight the following paper \cite{beznosikov2021decentralized}.
This work looks at a more general time-varying setting and suggests a new method. The upper bounds for their method are worse than for our method. We also mention papers on related topics:

\textit{Decentralized saddle point problems.}

The next work is devoted to centralized and decentralized distributed saddle problems \cite{beznosikov2020local}. It carries out lower bounds and optimal algorithms in the case when the communication network is constant (non-time-varying). See Table 1 for comparison our results for time-varying topology and results from \cite{beznosikov2020local} for constant network.

Also note the following works devoted to decentralized min-max problems. In \cite{liu2019decentralizedprox,rogozin2021decentralized} one can find algorithms for saddle point problems on fixed network. In paper \cite{rogozin2021decentralized}, it is near-optimal. Lower and upper bounds for decentralized min-max problems under data similarity condition are given in \cite{beznosikov2021distributed}. \cite{liu2019decentralized} studies the convergence of a decentralized methods for stochastic saddle point problems with homogeneous data on devices (all local functions $f_m$ are the same).

\textit{Minimization on time-varying networks.}

Decentralized methods are built upon combining iterations of classical first-order methods with communication steps. In the case of time-varying networks, a non-accelerated communication procedure is employed. Paper \cite{nedic2009distributed} can be named as an initial work on decentralized sub-gradient methods, and \cite{nedic2017achieving} proposed DIGing -- the first first-order minimization algorithm with linear convergence over time-varying networks. After that, PANDA, which is a dual method capable of working over time-varying graphs, was proposed in \cite{maros2018}. Analysis of DIGing and PANDA assumes that the underlying network is B-connected, that is, the union of B consequent networks is connected, while the network is allowed to be disconnected at some steps. Considering the time-varying graphs which stay connected at each iteration, decentralized Nesterov method \cite{rogozin2019optimal} has an accelerated rate under the condition that graph changes happen rarely enough, ADOM \cite{kovalev2021adom} and ADOM+ \cite{kovalev2021lower} are first-order optimization methods which achieve lower complexity bounds \cite{kovalev2021lower}. APM-C \cite{rogozin2019projected}, Acc-GT \cite{li2021accelerated} are accelerated methods over time-varying graphs, as well. The mentioned results are devoted to minimization algorithms and can be generalized to saddle-point problems. In this paper we generalize lower bounds of \cite{kovalev2021lower} to min-max problems and obtain an algorithm which reaches them up to a logarithmic factor.

\renewcommand{\arraystretch}{1.5}
\renewcommand{\tabcolsep}{10pt} 
\begin{table}[h!]
\vspace{-0.3cm}
\begin{center}
\begin{tabular}{ccc}
\hline
\multicolumn{1}{c}{} & \multicolumn{1}{c}{\textbf{time-varying network}} & \multicolumn{1}{c}{\textbf{constant network \cite{beznosikov2020local}}} \\ \hline
\multicolumn{1}{c}{}   & \multicolumn{2}{c}{{\tt lower}} \\ \hline
{\tt sc} & {$\Omega\left( R_0^2 \exp\left( - \frac{\mu K}{256L \chi} \right)\right)$} & {$\Omega\left( R_0^2 \exp\left( - \frac{\mu K}{128L \sqrt{\chi}} \right)\right)$} \\ \hline
{\tt c} & {$\Omega\left(\frac{L D^2 \chi}{K}\right)$} & {$\Omega\left(\frac{L D^2 \sqrt{\chi}}{K}\right)$} \\ \hline
\multicolumn{1}{c}{}   & \multicolumn{2}{c}{{\tt upper}} \\ \hline
{\tt sc} & {$\mathcal{\tilde O}\left( R_0^2 \exp\left( - \frac{\mu K}{8L \chi} \right)\right)$} & {$\mathcal{\tilde O}\left( R_0^2 \exp\left( - \frac{\mu K}{8L \sqrt{\chi}} \right)\right)$} \\ \hline
{\tt c} & {$\mathcal{\tilde O} \left(\frac{L D^2 \chi}{K}\right)$} & {$\mathcal{\tilde O}\left(\frac{L D^2 \sqrt{\chi}}{K} \right)$} \\ \hline
\end{tabular}
\vspace{0.4cm}
\caption{Lower and upper bounds for distributed smooth stochastic strongly-convex--strongly-concave ({\tt sc}) or convex-concave ({\tt c}) saddle-point problems in centralized and decentralized cases. Notation:
$L$ -- smothness constant of $f$, $\mu$ -- strongly-convex-strongly-concave constant, $R^2_0 = \|x_0 - x^* \|^2_2 + \|y_0 - y^* \|^2_2$, $D$ -- diameter of optimization set, $\chi$ -- condition number of communication graph (in time-varying case maximum of all graphs), $K$ -- number of communication rounds. In the case of upper bounds in the convex-concave case, the convergence is in terms of the "saddle-point residual", in the rest -- in terms of the (squared) distance to the solution.}
\end{center}
\label{table1}
\end{table}

\section{Preliminaries}

We use $\langle z,u \rangle := \sum_{i=1}^d z_i u_i$ to denote standard inner product of $z,u\in\R^d$. It induces $\ell_2$-norm in $\R^d$ in the following way $\|z\| := \sqrt{\langle z, z \rangle}$. We also introduce the following notation $\text{proj}_{\mathcal{Z}}(z) = \min_{u \in \mathcal{Z}}\| u - z\|$ -- the Euclidean projection onto $\mathcal{Z}$.

We work with the problem \eqref{SPP}, where the sets $\mathcal{X} \subseteq \mathbb{R}^{n_x}$ and $\mathcal{Y} \subseteq \mathbb{R}^{n_y}$ are convex sets. Additionally, we introduce the set $\mathcal{Z} = \mathcal{X} \times \mathcal{Y}$, $z = (x,y)$ and the operator $F$:
\begin{equation}
\label{opSP}
    F_m(z) = F_m(x,y) = \begin{pmatrix}
\nabla_x f_m(x,y)\\
-\nabla_y f_m(x,y)
\end{pmatrix}.
\end{equation}
This notation is needed for shortness. 

\textbf{Problem setting.} Next, we introduce the following assumptions:

\textit{Assumption 1(g).}
 $f(x,y)$ is $L$ - smooth, if for all $z_1, z_2 \in \mathcal{Z}$
    \begin{eqnarray}\label{as1g}
    \|F(z_1) - F(z_2)\| \leq L\|z_1-z_2\|.\end{eqnarray}
    
\textit{Assumption 1(l).}
 For all $m$, $f_m(x,y)$ is Lipschitz continuous with constant $L_{\max}$, it holds that for all $z_1, z_2 \in \mathcal{Z}$
    \begin{eqnarray}\label{as1l}
    \|F_m(z_1) - F_m(z_2)\| \leq L_{\max}\|z_1-z_2\|.\end{eqnarray}
    
\textit{Assumption 2(s).}
$f(x,y)$ is strongly-convex-strongly-concave with constant $\mu$, if for all $z_1, z_2 \in \mathcal{Z}$
    \begin{eqnarray} \label{as2g}\langle F(z_1) - F(z_2), z_1 - z_2 \rangle \geq \mu\|z_1-z_2\|^2.\end{eqnarray}
    
\textit{Assumption 2(c).}
$f(x,y)$ is convex-concave, if $f(x,y)$ is strongly-convex-strongly-concave with $0$.

\textit{Assumption 3.} $\mathcal{Z}$ -- compact bounded, i.e. for all $z, z'\in \mathcal{Z}$
    \begin{eqnarray} \label{as5}
    \| z - z'\| \leq D.
    \end{eqnarray}
All assumptions are standard in the literature.

\textbf{Network setting.} In each moment of time (iteration) $t$, the communication network  is modeled  as a connected, undirected graph graph $\mathcal{G}(t) \triangleq (\mathcal{V}, \mathcal{E}(t))$, where $\mathcal{V} := \{1,\ldots, M\}$ denotes  the vertex    set--the set of devices (does not change in time) and  $\mathcal{E}(t) := \{(i,j) \, |\, i,j \in \mathcal{V} \}$ represents the set of edges--the communication links at the moment $t$;  $(i,j) \in \mathcal{E} (t)$ iff there exists a communication link between devices $i$ and $j$ in moment $t$.  

As mentioned earlier, the gossip protocol is the most popular communication procedures in decentralized
setting. This approach uses a certain matrix $W$. Local vectors during communications are "weighted" by multiplication of a vector with $W$. The convergence of decentralized algorithms is determined by the properties of this matrix. Therefore, we introduce the following assumption:

\textit{Assumption 4.} We call a matrix $W(t)$ a gossip matrix at the moment $t$ if it satisfies the following conditions: 1) $W(t)$ is an $M \times M$ symmetric, 2) $W(t)$ is positive semi-definite, 3) the kernel of $W(t)$ is the set of constant vectors, 4) $W(t)$ is defined on the edges of the network at the moment $t$: $W_{ij}(t) \neq 0$ only if $i=j$ or $(i,j) \in \mathcal{E}(t)$.

Let $\lambda_1(W(t)) \geq \ldots \geq \lambda_M(W(t)) = 0$ be the spectrum of $W(t)$, and define condition number $\chi = \max_t \chi(W(t)) = \max_t \frac{\lambda_1(W(t))}{\lambda_{M-1}(W(t))}$. Note that in practice we use not the matrix $W(t)$, but $\tilde W(t) = I - \frac{W(t)}{\lambda_1(W(t))}$, since this type of matrices are used in consensus algorithms \cite{boyd2006randomized}. To estimate the convergence speed, we introduce 
\begin{align*}
    \rho &= \max_t \lambda_2(\tilde W(t)) = \max_t \left[1 - \frac{\lambda_{M-1}(W(t))}{\lambda_1(W(t))}\right] = \max_t \left[1 - \frac{1}{\chi(W(t))}\right] \\
    &= 1 - \frac{1}{\max_t \chi(W(t))} = 1 - \frac{1}{\chi}.
\end{align*}

\section{Main part}

We divide our contribution into two main parts, first we discuss lower bounds for decentralized saddle point problems over time-varying graphs. In the second part, we present an algorithm that achieves the lower bounds (up to logarithmic factors and numerical constants).

\subsection{Lower bounds}

Before presenting lower bounds, we must restrict the class of algorithms for which our lower bounds are valid. For this we introduce the following black-box procedure.

\begin{definition} \label{proc}
    Each device $m$ has its own local memories $\mathcal{M}^x_{m}$ and $\mathcal{M}^y_{m}$ for the $x$- and $y$-variables, respectively--with  initialization $\mathcal{M}_{m}^x = \mathcal{M}_{m}^y= \{0\}$.    $\mathcal{M}_{m}^x$ and $\mathcal{M}_{m}^x$ are updated as follows:\\
   
    $\bullet$ \textbf{Local computation:} Each device $m$ computes and adds to its $\mathcal{M}^x_{m}$ and $\mathcal{M}^y_{m}$ a finite number of points $x,y$, each satisfying 
    \begin{equation}\begin{aligned}\label{eq:oracle-opt-step}
        x \in \text{span} \big\{x'~,~\nabla_x f_m(x'',y'')\big\},\quad 
        y \in \text{span} \big\{y'~,~\nabla_y f_m(x'',y'')\big\}, 
    \end{aligned}\end{equation}
    for given $x', x'' \in \mathcal{M}^x_{m}$ and  $y', y'' \in \mathcal{M}^y_{m}$.

    $\bullet$ \textbf{Communication:} Based upon communication round among neighbouring nodes at the moment $t$,   $\mathcal{M}^x_{m}$ and $\mathcal{M}^y_{m}$ are updated according to
    \begin{equation}\label{eq:oracle-comm}
        \mathcal{M}^x_{m} := \text{span}\left\{\bigcup_{(i,m) \in \mathcal{E}(t)} \mathcal{M}^x_{i} \right\}, \quad 
        \mathcal{M}^{y}_{m} := \text{span}\left\{\bigcup_{(i,m) \in \mathcal{E}(t)} \mathcal{M}^y_{i} \right\}.
    \end{equation}

    $\bullet$ \textbf{Output:} 
    The final global output at the current moment of time is calculated as: 
    \begin{align*}
        x \in \text{span}\left\{\bigcup_{m=1}^M \mathcal{M}^x_{m} \right\},~~y \in \text{span}\left\{\bigcup_{m=1}^M \mathcal{M}^y_{m} \right\}.
    \end{align*}
\end{definition} 

This definition includes all algorithms capable of making local gradient updates, as well as exchanging information with neighbors.
Notice that the proposed oracle builds on   \cite{scaman2017optimal}  for minimization problems over networks. 

\begin{theorem}\label{th-LB-distributed}
For any $L$ and $\mu$ , there exists a  saddle point problem in the form (\ref{SPP}) with $\mathcal{Z}=\mathcal{R}^{2d}$(where $d$ is sufficiently
large) and non-zero solution $y^*$. All local functions $f_m$ of this problem are  $L$-smooth, $\mu$-strongly-convex-strongly-concave. Then, for any $\chi \geq 1$, there exists a sequence of gossip matrices $W(t)$ over the connected (at each moment) graph  $\mathcal{G} (t)$,  satisfying Assumption 4 with condition number $\chi$, such that any decentralized algorithm satisfying Definition~\ref{proc} and using the gossip matrices $W(t)$  produces the following estimate on the global output $z=(x,y)$ after $K$ communication rounds:
\begin{equation*}
    \|z^{K} - z^*\|^2 = \Omega\left(\exp\left( - \frac{256\mu}{L - \mu} \cdot \frac{K}{\chi} \right) \| y^*\|^2\right).
\end{equation*}
\end{theorem}

The idea of finding lower bounds is to construct an example of "bad" functions and the "critical" location of these functions on the nodes. In papers \cite{beznosikov2020local,beznosikov2021distributed}, lower bounds for decentralized saddle point problems (but on fixed communication networks) were already investigated. Examples of "bad" functions and their analysis can be taken from these works. An example of "bad" time-varying topology of the node connection is a star with a changing center. Obtaining lower bounds using such varying networks for minimization problems was obtained in \cite{kovalev2021lower}. To prove Theorem 1 we need to combine results \cite{beznosikov2020local} and \cite{kovalev2021lower}.

The following statement interprets Theorem 1 in terms of the number of local computations on each device and the number of communications between them.
\begin{corollary}\label{cor1}
In the setting of Theorem~\ref{th-LB-distributed}, the number of communication rounds required to obtain a $\varepsilon$-solution is lower bounded by\vspace{-0.2cm}
\begin{equation*}
    \Omega\left( \chi \frac{L}{\mu} \cdot  \log \left(\frac{\| y^*\|^2}{\varepsilon}\right)\right).
\end{equation*}
Additionally, we can get a lower bound for the number of local calculations on each of the devices:
\begin{equation*}
    \Omega\left( \frac{L}{\mu} \cdot  \log \left(\frac{\| y^*\|^2}{\varepsilon}\right)\right).
\end{equation*}
\end{corollary}

Also we want to find lower bounds for the case of  (non strongly) convex-concave problems, one can use  regularization and consider the following objective function 
\begin{align*}
    f(x,y) + \frac{\varepsilon}{4D^2} \cdot \|x - x^0 \|^2 - \frac{\varepsilon}{4D^2}\cdot \|y - y^0 \|^2,
\end{align*}
which is   strongly-convex-strongly-concave with constant  $\mu = \frac{\varepsilon}{2D^2}$, where $\varepsilon$ is a precision of the solution and $D$ is the  diameter of the  sets $\mathcal{X}$ and $\mathcal{Y}$. The resulting new SPP problem is solved to $\varepsilon/2$-precision in order to guarantee an accuracy $\varepsilon$ in computing the solution of  the original problem. Therefore, we can easily deduce the lower bounds for convex-concave case
\begin{equation*}
\label{lower_cc}
    \Omega\left( \chi \frac{L D^2}{\varepsilon} \right)~\text{communication rounds}\quad \text{and}\quad  \Omega\left(\frac{L D^2}{\varepsilon}\right) ~\text{local computations}.
\end{equation*}
See Table 1 to compare with lower bounds for constant networks. 

\subsection{Near-optimal algorithm}

In this part, we present an Algorithm that achieves lower bounds (up to logarithmic terms). Our Algorithm uses an auxiliary procedure for communication. This is a classic procedure - Gossip Algorithm.

\begin{algorithm} [th]
	\caption{Gossip Algorithm ({\tt Gossip})}
	\label{alg3}
	\begin{algorithmic}
\State
\noindent {\bf Parameters:} Vectors $z_1, ..., z_M$, communic. rounds $H$.
\State \noindent {\bf Initialization:}
Construct matrix $\textbf{z}$ with rows $z^T_1, ..., z^T_M$.
\State Choose $\textbf{z}^0 = \textbf{z}$.
\For {$h=0,1, 2, \ldots, H$ } 
\State $\textbf{z}^{h+1} = \tilde W (h) \cdot  \textbf{z}^{h}$
\EndFor
\State
\noindent {\bf Output:} rows $z_1, ..., z_M$ of  $\textbf{z}^{H+1}$ .
	\end{algorithmic}
\end{algorithm}

The essence of the {\tt Gossip} is very simple. Initially, there are vectors $z_1$ and $z_M$, which are stored on their devices. Our goal is to get a vector close to the $\bar z = \frac{1}{M} \sum_{m=1}^M z_m$ vector on all devices. At each iteration, each device exchange local vectors with its neighbors, and then modify its local vector by averaging local vector and vectors of neighbors with weights from the matrix $W (h)$.

We are now ready to present our main algorithm. It is based on the classical method for smooth saddle point problems - Extra Step Method (Mirror Prox) \cite{Nemirovski2004,juditsky2008solving}. With the right choice of $H$, we can achieve averaging of all vectors with good accuracy. In particular, we can assume that $z^k_1 \approx \ldots \approx z^k_M$. For more details about the choice of $H$ and a detailed analysis of the algorithm (taking into account that in the general $z^k_1 \neq \ldots \neq z^k_M$), see in the full version of the paper.

\begin{algorithm} [th]
	\caption{Time-Varying Decentralized Extra Step Method ({\tt TVDESM})}
	\label{alg2}
	\begin{algorithmic}
\State
\noindent {\bf Parameters:}  Stepsize $\gamma \leq \frac{1}{4L}$, number of {\tt Gossip} steps $H$.\\
\noindent {\bf Initialization:} Choose  $(x^0,y^0)=z^0\in \mathcal{Z}$, $z^0_m = z^0$.
\For {$k=0,1, 2, \ldots, $ } 
 \State Each machine $m$ computes  \hspace{0.1cm} $\hat z_m^{k+1/2} = z_m^{k} - \gamma \cdot F_m(z^k_m)$
 \State  Communication: \hspace{0.1cm} $\tilde z^{k+1/2}_1, \ldots,  \tilde z^{k+1/2}_M$ ={\tt Gossip}$(\hat z^{k+1/2}_1, \ldots, \hat z^{k+1/2}_M, H)$
\State Each machine $m$ computes $z^{k+1/2}_m = \text{proj}_{\mathcal{Z}}(\tilde z^{k+1/2}_m)$,
\State Each machine $m$ computes \hspace{0.1cm} $\hat z_m^{k+1} = z_m^{k} - \gamma \cdot F_m(z^{k+1/2}_m)$
\State Communication: \hspace{0.1cm} $\tilde z^{k+1}_1, \ldots,  \tilde z^{k+1}_M$ ={\tt Gossip}$(\hat z^{k+1}_1, \ldots, \hat z^{k+1}_M, H)$
\State Each machine $m$ computes \hspace{0.1cm} $z^{k+1}_m = \text{proj}_{\mathcal{Z}}(\tilde z^{k+1}_m)$
\EndFor
	\end{algorithmic}
\end{algorithm}

The analysis of Algorithm \ref{alg2} is derived from the analysis of classical extrastep method. We study the convergence properties of sequence $\{\bar z^k\}_{k=0}^\infty$, where $\bar{z}^k = \frac{1}{M} \sum_{m=1}^M z_m^k$. Note that $\bar z^k$ is not held at any agent; instead, this quantity is only used in the analysis. Algorithm \ref{alg2} employs gossip averaging after each extra-step. Therefore, the method uses an approximate value of $F(\bar z^k)$ when performing updates, and the approximation error is driven by the number of gossip iterations $H$. The analysis of Algorithm \ref{alg2} comes down to studying extrastep method which uses inexact values of $F$ at each iteration. Given a target accuracy $\varepsilon$, we choose the number of gossip iterations $H$ proportional to $\varepsilon$. Since \texttt{Gossip} (Algorithm \ref{alg3}) is a linearly convergent method, $H$ is proportional to $\log(1/\varepsilon)$. As a result, we have a $\log^2(1/\varepsilon)$ term in the number of communication rounds of Algorithm \ref{alg2}.

\begin{theorem} 
Let $\{ z_m^k\}^K_{k \geq 0}$ denote the iterates of Algorithm~\ref{alg2} for solving problem \eqref{SPP} after $K$ communication rounds. Let Assumptions 1(g,l) and 4 be satisfied. Then, if $\gamma \leq \frac{1}{4L}$, we have the following estimates in

$\bullet$ $\mu$-strongly-convex--strongly-concave case (Assumption 2(s)):
\begin{eqnarray*}
 \| \bar z^{K+1}\!-\!z^* \|^2\!=\!\mathcal{\tilde O}\left( \| z^{0} - z^* \|^2 \exp\left( -\frac{\mu K}{8L{\chi}} \right) \right) ,
 \end{eqnarray*}

$\bullet$ convex--concave case (Assumption 2 and 3):
\begin{equation*}
    \text{gap}(\bar z^{K + 1}_{avg}) = \mathcal{\tilde O}\left(\frac{L \Omega_z^2 \chi}{K} \right),
\end{equation*}

where $\bar z^{t} = \frac{1}{M} \sum\limits_{m=1}^M z_m^{t}$, $\bar z^{k+1}_{avg} = \frac{1}{M(k+1)} \sum\limits_{t=0}^k \sum\limits_{m=1}^M z_m^{t+1/2}$ and $$\text{gap}(z) = \max_{y'\in \mathcal{Y}} f(x, y') -  \min_{x'\in \mathcal{X}} f(x', y).$$
\end{theorem}

\begin{corollary}\label{cor2}
In the setting of Theorem~2, the number of communication rounds required for Algorithm 2 to obtain a $\varepsilon$-solution is upper bounded by\vspace{-0.2cm}
\begin{equation*}
    \mathcal{\tilde O}\left( \chi \frac{L}{\mu}\right)
\end{equation*}
in $\mu$-strongly-convex--strongly-concave case and 
\begin{equation*}
\mathcal{\tilde O}\left( \chi \frac{L D^2}{\varepsilon} \right)
\end{equation*}
in convex-concave case.
Additionally, one can obtain upper bounds for the number of local calculations on each of the devices:
\begin{equation*}
    \mathcal{O}\left( \frac{L}{\mu} \cdot  \log \left(\frac{\| z^0 - z^*\|^2}{\varepsilon}\right)\right)
\end{equation*}
in $\mu$-strongly-convex--strongly-concave case and 
\begin{equation*}
\mathcal{O}\left(\frac{L D^2}{\varepsilon} \right)
\end{equation*}
in convex-concave case.
\end{corollary}

Corollary \ref{cor2} illustrates that Algorithm \ref{alg2} achieves lower bounds both for convex-concave and $\mu$-strongly-convex-strongly-concave cases up to a logarithmic factor (the lower bounds are determined in Corollary \ref{cor1}). It can be observed that the complexity bounds are constituted of two factors: $\chi$ representing network connectivity and $L/\mu$ of $LD^2/\varepsilon$ corresponding to the objective function. This effect is typical for decentralized optimization (see i.e. \cite{scaman2017optimal}). On the contrary to distributed minimization tasks, the dependence on function condition number $L/\mu$ is unimprovable for min-max problems (i.e. this factor cannot be enhanced to $\sqrt{L/\mu}$). Moreover, the dependence on $\chi$ cannot be improved to $\sqrt\chi$, since we focus on time-varying networks \cite{kovalev2021lower}.

\section{Conclusion}

In conclusion, we briefly summarize the contributions of this paper and discuss the directions for future work. Our findings consist of two parts: lower bounds and optimal (up to a logarithmic factor) algorithms.

First, we derived the lower bounds for the classes of convex-concave and strongly-convex-strongly-concave min-max problems over time-varying graphs. The graph is assumed to be connected at each communication round. However, we studied only one class of time-varying networks. Other classes are connected to different assumptions on the network structure. In particular, in B-connected networks \cite{nedic2017achieving} the graph can be disconnected at some times, but the union of any B consequent graphs must be connected. Yet another possible assumption is the randomly changing graph with a contraction property of $W$ in expectation \cite{koloskova2020unified}. Developing lower bounds for min-max problems for these two classes is an open question in decentralized optimization.

Second, we proposed a near-optimal algorithm with a gossip subroutine resulting in squared logarithmic factor. Developing an algorithm without an additional logarithmic factor would close the gap in theory and result in a more practical algorithm with less parameters to fine-tune. Possible directions for developing such an algorithm are generalizations of dual-based approaches for minimization \cite{kovalev2021adom,maros2018} and gradient-tracking \cite{nedic2017achieving,maros2018}.

Finally, the comparison of our algorithm to existing works requires additional numerical experiments, which is left for future work.

\bibliographystyle{splncs04}
\bibliography{ltr}

\end{document}